\documentclass[12pt]{amsart}
\usepackage{amsmath,amsthm,amsfonts,amssymb}
\usepackage[
margin=2.5cm]{geometry}
\usepackage[usenames,dvipsnames]{color}
\usepackage{mathrsfs}
\usepackage[colorlinks,allcolors=blue]{hyperref}
\usepackage{tikz,tikz-qtree}
\usetikzlibrary{matrix,arrows,trees,shapes}
\usepackage{xcolor}
\usepackage{enumitem}
\setlist[enumerate]{%
  leftmargin=*,%
  label=(\roman*)} 
\setlist[itemize]{leftmargin=*}
\usepackage[mathscr]{eucal}
\usepackage{multirow, array}

\usepackage[yyyymmdd]{datetime}

\usepackage{fancyhdr}
\newtheorem{theorem}{Theorem}[section]
\newtheorem{lemma}[theorem]{Lemma}
\newtheorem{proposition}[theorem]{Proposition}
\newtheorem{corollary}[theorem]{Corollary}
\newtheorem{propdef}[theorem]{Proposition/Definition}

\theoremstyle{definition}%
\newtheorem{example}[theorem]{Example}
\newtheorem{remark}[theorem]{Remark}
\newtheorem{algorithm}[theorem]{Algorithm}
\newtheorem{notation}[theorem]{Notation}

\DeclareMathOperator{\hilb}{Hilb}
\DeclareMathOperator{\PP}{\mathbb{P}}
\DeclareMathOperator{\hilbtree}{\mathscr{H}}

\newcommand{\hf}{\mathsf{h}}
\newcommand{\hp}{\mathsf{p}}
\newcommand{\hq}{\mathsf{q}}
\newcommand{\lift}{\Phi}
\newcommand{\plus}{A}

\newcommand{\kk}{\mathbb{K}}
\newcommand{\NN}{\mathbb{N}}
\newcommand{\QQ}{\mathbb{Q}}
\newcommand{\ZZ}{\mathbb{Z}}
\newcommand{\FF}{\mathbb{F}}

\DeclareMathOperator{\lexg}{\succ}

\newcommand{\relphantom}[1]{\mathrel{\phantom{ #1 }}}
\newcommand{\llrr}[1]{ \langle #1 \rangle } 

\begin{document}

\begin{abstract}
  We extend the recent classification of Hilbert schemes with two
  Borel-fixed points to arbitrary characteristic.  We accomplish this
  by synthesizing Reeves' algorithm for generating strongly stable
  ideals with the basic properties of Borel-fixed ideals and our
  previous work classifying Hilbert schemes with unique Borel-fixed
  points.
\end{abstract}

\title{Hilbert Schemes with Two Borel-fixed Points in Arbitrary
  Characteristic}
\author{Andrew~P.~Staal}
\email{andrew.staal@uwaterloo.ca}
\address{ Department of Pure Mathematics \\
  University of Waterloo \\
  200 University Avenue West \\
  Waterloo, Ontario, Canada, N2L 3G1 }%
\date{}
\maketitle

\section{Introduction}

Hilbert schemes parametrizing closed subschemes with a fixed Hilbert
polynomial in projective space are fundamental moduli spaces.  In
\cite{Staal--2020}, we provide a rough guide to aid in predicting the
complexity of the geometry of Hilbert schemes.  Specifically, we
classify Hilbert schemes with unique Borel-fixed points, over an
algebraically closed or characteristic $0$ field, and show that these
Hilbert schemes are smooth and irreducible (they are also rational).
Recently, Ramkumar \cite{Ramkumar--2019} has built on this,
classifying Hilbert schemes with two Borel-fixed points and analyzing
their local geometry, when the base field has characteristic $0$.
Here we show that the classification holds (with one
minor adjustment) over arbitrary
infinite fields.

Let $\hilb^{\hp}(\PP^n)$ be the Hilbert scheme parametrizing closed
subschemes of $\PP^n_{\kk}$ with Hilbert polynomial $\hp$, where $\kk$
is an
infinite field.  Macaulay classified Hilbert polynomials of
homogeneous ideals in \cite{Macaulay--1927}.  Any such admissible
Hilbert polynomial $\hp(t)$ has a unique combinatorial expression of
the form $\sum_{j=1}^r \binom{t+b_j-j+1}{b_j}$, for integers $b_1 \ge
b_2 \ge \dotsb \ge b_r \ge 0$, with $d := \deg \hp = b_1$.  Our main
result is the following.

\begin{theorem}
  \label{thm:main1}
  Let $n \ge d + 2$.  The Hilbert scheme $\hilb^{\hp}(\PP^n)$ has
  exactly two Borel-fixed points if and only if one of the following
  conditions holds:
  \begin{enumerate}
  \item
    \begin{enumerate}[label={(\alph*)}]
    \item $b_1 = 0$ and $r = 3$,

    \item[(a')] $b_1 = 0$ and $r = 4$, if $n = 2$ and $\kk$ does not
      have characteristic $2$,
    
    \item $b_1 = \dotsb = b_{r-1} = 1$ and $b_r = 0$, for $r-1 \neq 1, 3$,

    \item $b_1 = \dotsb = b_{r-1} \ge 2$ and $b_r = 0$, for $r-1 \neq
      1$,
    \end{enumerate}

  \item
    \begin{enumerate}[label={(\alph*)}]
    \item $b_1 > b_2 = 0$ and $r = 3$,
      
    \item $b_1 = \dotsb = b_{r-2} > b_{r-1} = 1$ and $b_r = 0$, for
      $r-2 \neq 2$, and

    \item $b_1 = \dotsb = b_{r-2} > b_{r-1} \ge 2$ and $b_r = 0$.
    \end{enumerate}
  \end{enumerate}
\end{theorem}

A point on $\hilb^{\hp}(\PP^n)$ is Borel-fixed if its corresponding
(saturated) ideal is fixed by the linear action of the Borel group of
upper-triangular matrices in $\operatorname{GL}_{n+1}(\kk)$ on the
polynomial ring $\kk[x_0, x_1, \dotsc, x_n]$.  Borel-fixed ideals play
an important geometric role, because they often function as markers
for interesting local geometry on Hilbert schemes.  Many fundamental
properties of Hilbert schemes have been understood by passing from a
homogeneous ideal to its generic initial ideal, which is Borel-fixed
\cite{Bayer--Stillman--1987}.  Thus, to analyze a chosen Hilbert
scheme, it is beneficial to understand all of its Borel-fixed points.
In characteristic $0$, Borel-fixed ideals are strongly stable.  The
combinatorial criterion defining strongly stable ideals is robust
enough that many of their properties are well understood.  For
example, if $I$ is strongly stable, then the minimal free resolution
of $I$ is an iterated mapping cone \cite{Eliahou--Kervaire--1990,
  Peeva--2011}, the regularity of $I$ is the largest of the degrees of
its minimal monomial generators \cite{Bayer--Stillman--1987}, and the
saturation of $I$ is generated by Reeves' algorithm
\cite{Reeves--1992, Moore--Nagel--2014,
  Cioffi--Lella--Marinari--Roggero--2011}.  This makes studying
Borel-fixed ideals much easier than in general.

Strongly stable ideals, including lexicographic ideals, are always
Borel-fixed---their definition is derived by describing Borel-fixed
ideals in characteristic $0$.  However, in positive characteristic,
Borel-fixed ideals generally satisfy a looser
combinatorial condition and are harder to understand.  Progress has
been made \cite{Pardue--1994, Herzog--Popescu--2001,
  Sinefakopoulos--2007}, but the minimal free resolution of a
\emph{nonstandard} Borel-fixed ideal $I$ (see \S \ref{ch:borelfixed})
is only understood in limited cases, while obtaining a nonminimal free
resolution or the regularity of $I$ can require finding a Pommaret
basis of $I$ \cite{Seiler--2009,
  Albert--Fetzer--Saenz-de-Cabezon--Seiler--2015}.  An algorithm to
generate Borel-fixed ideals via Pommaret bases appears in
\cite{Bertone--2015} (see Remark \ref{rmk:Bertone}).

Theorem~\ref{thm:main1} provides a class of Hilbert schemes with the
notable feature that their Borel-fixed points are all strongly stable
(as does the main theorem in \cite{Staal--2020}).  To obtain this
classification, we rely only on Pardue's description of the elementary
properties of Borel-fixed ideals over infinite fields
\cite{Pardue--1994}, in conjunction with our classification in
\cite{Staal--2020}.  Despite the absence of
the theory of expansions \cite{Moore--Nagel--2014}
for Borel-fixed ideals in positive characteristic, basic recursive
aspects of Reeves' algorithm are sufficient to generate the
Borel-fixed ideals on the Hilbert schemes we consider.  In particular,
once we understand the strongly stable points on our Hilbert schemes
and the Borel-fixed points on certain ``nearby'' ones, the weaker
exchange property of Borel-fixed ideals is sufficient to complete our
classification.

Note that when $\kk$ has characteristic $0$, the local deformation
theory at these Borel-fixed ideals is worked out in
\cite{Ramkumar--2019}.  The Hilbert schemes described in Theorem
\ref{thm:main1}(i) also appear in the recent classification of smooth
Hilbert schemes \cite{Skjelnes--Smith--2020} (see Remark
\ref{rmk:SkSm}).

\subsection*{Conventions.}

Throughout, $\kk$ denotes an
infinite field, $\NN$ denotes the nonnegative integers, $S := \kk[x_0,
  x_1, \dotsc, x_n]$ is the standard $\ZZ$-graded polynomial ring, and
$\mathfrak{m}_k := \langle x_0, x_1, \dotsc, x_k \rangle$ is its
homogeneous ideal, for $0 \le k \le n$.  The Hilbert function and
polynomial of the quotient $S / I$ by a homogeneous ideal $I$ are
denoted $\hf_I$ and $\hp_I$, respectively.

\subsection*{Acknowledgments}

We thank Joachim Jelisiejew for pointing out the overlap between
Ritvik Ramkumar's work and our own, and Greg Smith for helpful
suggestions.  We thank Eran Nevo and The Hebrew University of
Jerusalem for supporting this research, as well as the Topology and
Geometry Group at the University of Waterloo.

\section{Background}

In this section, we provide a concise exposition of some key ideas
from \cite{Staal--2020} needed to address the ``two Borels'' case.
Particularly relevant is the arrangement of admissible Hilbert
polynomials by their associated partitions into an infinite full
binary tree, and the induced arrangement of Hilbert schemes into
copies of this tree indexed by the codimensions of the parametrized
schemes.

\subsection{Admissible Hilbert Polynomials}

Let $S := \kk[x_0, x_1, \dotsc, x_n]$ denote the homogeneous
coordinate ring of $\PP^n_{\kk}$ over an
infinite field $\kk$.  Let $M$ be a finitely generated graded
$S$\nobreakdash-module.  The \emph{\bfseries Hilbert function} $\hf_M
\colon \ZZ \to \ZZ$ of $M$ is defined by $\hf_M(i) := \dim_{\kk} M_i$,
for all $i \in \ZZ$.  Every such $M$ has a \emph{\bfseries Hilbert
polynomial} $\hp_M(t) \in \QQ[t]$ such that $\hf_M(i) = \hp_M(i)$, for
$i \gg 0$; see \cite[Theorem~4.1.3]{Bruns--Herzog--1993}.  For a
homogeneous ideal $I \subset S$, we denote by $\hf_I$ and $\hp_I$ the
Hilbert function and Hilbert polynomial \emph{of the quotient module}
$S/I$, respectively.

\begin{notation}
  \label{rmk:binconv}
  For integers $j,k$, set $\binom{j}{k} = \frac{j!}{k!(j-k)!}$ if $j
  \ge k \ge 0$ and $\binom{j}{k} = 0$ otherwise.  For a variable $t$
  and $a, b \in \ZZ$, define $\binom{t + a}{b} = \frac{(t + a)(t + a -
    1) \dotsb (t + a - b + 1)}{b!} \in \QQ[t]$ if $b \ge 0$, and
  $\binom{t + a}{b} = 0$ otherwise.
\end{notation}

A polynomial is an \emph{\bfseries admissible Hilbert polynomial} if
it is the Hilbert polynomial of a closed subscheme in some $\PP^n$.
Admissible Hilbert polynomials correspond to nonempty Hilbert schemes;
we always work with nonzero admissible Hilbert polynomials.  We use
the following well-known classification due to Macaulay.

\begin{proposition} 
  \label{prop:expressions}
  The following conditions are equivalent:
  \begin{enumerate}
  \item The polynomial $\hp(t) \in \QQ[t]$ is a nonzero admissible
    Hilbert polynomial.
    
  \item There exist integers $e_0 \ge e_1 \ge \dotsb \ge e_d > 0$ such
    that $\hp(t) = \sum_{i=0}^{d} \binom{t + i}{i + 1} - \binom{t + i
      - e_i}{i + 1}$.

  \item There exist integers $b_1 \ge b_2 \ge \dotsb \ge b_r \ge 0$
    such that $\hp(t) = \sum_{j=1}^r \binom{t + b_j - j+1}{b_j}$.
  \end{enumerate}
  Moreover, these correspondences are bijective.
\end{proposition}

\begin{proof} $\;$
  \begin{enumerate}[leftmargin=2cm]
  \item[(i) $\Leftrightarrow$ (ii)] This is proved in
    \cite{Macaulay--1927}; see also \cite[Corollary~3.3 and
      Corollary~5.7]{Hartshorne--1966}.

  \item[(i) $\Leftrightarrow$ (iii)] This follows from
    \cite[Erinnerung~2.4]{Gotzmann--1978}; see also
    \cite[Exercise~4.2.17]{Bruns--Herzog--1993}.
  \end{enumerate}
  Uniqueness of the sequences of integers also follows.
\end{proof}

The \emph{\bfseries Macaulay expression} of an admissible Hilbert
polynomial $\hp$ is its unique expression
\[
\textstyle \hp(t) = \sum_{i=0}^{d} \binom{t + i}{i + 1} - \binom{t + i
  - e_i}{i + 1},
\]
where $e_0 \ge e_1 \ge \dotsb \ge e_d > 0$, and the \emph{\bfseries
Gotzmann expression} of $\hp$ is its unique expression
\[
\textstyle \hp(t) = \sum_{j=1}^{r} \binom{t + b_j - j+1}{b_j},
\]
for $b_1 \ge b_2 \ge \dotsb \ge b_r \ge 0$.  From these, we can read
off the degree $d = b_1$, the leading coefficient $e_d / d!$, and the
\emph{\bfseries Gotzmann number} $r$ of $\hp$; the latter bounds the
Castelnuovo--Mumford regularity of saturated ideals with Hilbert
polynomial $\hp$.  In particular, such ideals are generated in degrees
up to $r$; see \cite[p.\ 300--301]{Iarrobino--Kanev--1999}.  If
$\hp(t)$ has Macaulay and Gotzmann expressions as above, then $r =
e_0$ and the (nonnegative) partition $(b_1, b_2, \dotsc, b_r)$ is
conjugate to the partition $(e_1, e_2, \dotsc, e_d)$; see
\cite[Lemma~2.4]{Staal--2020}.  In particular, the partition $(b_1,
b_2, \dotsc, b_r)$ has $e_i - e_{i+1}$ parts equal to $i$, for all $0
\le i \le d$.  We call $(e_0, e_1, \dotsc, e_d)$ the \emph{\bfseries
Macaulay partition} of $\hp$ and $(b_1, b_2, \dotsc, b_r)$ the
\emph{\bfseries Gotzmann partition} of $\hp$.  We use Gotzmann
expressions of admissible Hilbert polynomials from here on.

There are two fundamental binary relations on admissible Hilbert
polynomials.  To describe these, let $\hp$ have Gotzmann partition
$(b_1, b_2, \dotsc, b_r)$.  The first ``lifts'' or ``integrates''
$\hp$ to the polynomial $\lift(\hp)$ with partition $(b_1 +1, b_2 +1,
\dotsc, b_r + 1)$.  The second ``adds one,'' taking $\hp$ to
$\plus(\hp) := 1 + \hp$, which has partition $(b_1, b_2, \dotsc, b_r,
0)$.  These roughly correspond to the geometric operations of forming
a cone over a subscheme of $\PP^{n}$ (with cone point outside $\PP^n$
in an ambient $\PP^{n+1}$) and forming the disjoint union with a
reduced point (inside $\PP^n$), respectively.  Both $\lift(\hp)$ and
$\plus(\hp)$ are admissible by Proposition~\ref{prop:expressions}.

\subsection{A Geography of Hilbert Schemes}
\label{ch:geo}

Using the binary relations $\lift$ and $\plus$, we observe in
\cite{Staal--2020} that the set of (nonempty) Hilbert schemes of the
form $\hilb^{\hp}(\PP^n)$ can be thought of as the vertex set of a
collection of infinite full binary trees.

\begin{propdef}
  \label{propdef:hilbtree}
  For each positive codimension $c$, let $\hilbtree_c$ be the graph
  whose vertices are all nonempty Hilbert schemes $\hilb^{\hp}(\PP^n)$
  parametrizing codimension $c$ subschemes and whose edges are all
  pairs $\bigl( \hilb^{\hp}(\PP^n), \hilb^{\plus(\hp)}(\PP^n) \bigr)$
  and $\bigl( \hilb^{\hp}(\PP^n), \hilb^{\lift(\hp)}(\PP^{n+1})
  \bigr)$, where $\hp$ is an admissible Hilbert polynomial and $n = c
  + \deg \hp$.  Then $\hilbtree_c$ is an infinite full binary tree
  with root $\hilb^1(\PP^c) = \PP^c$.
\end{propdef}

\begin{proof}
  See \cite[Theorem~2.14]{Staal--2020}.
\end{proof}

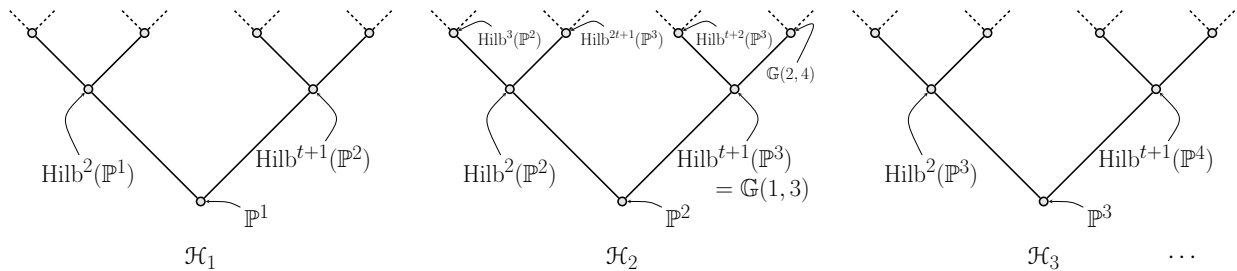
\begin{figure}[!ht] 
  \centering
  \resizebox{\textwidth}{!}{%
    \begin{tikzpicture}
      \tikzstyle{point}=[circle, ultra thick, draw=black, fill=gray!20,
        inner sep=3pt, minimum width=4pt, minimum height=4pt]

      \node (n1) at (0,-2) {\Huge $\hilbtree_1$};
      \node (o)[point] at (0,0) {};
      \node (a)[point] at (-4,4) {};
      \node (b)[point] at (4,4) {};
      \node (aa)[point] at (-6,6) {};    
      \node (ab)[point] at (-2,6) {};
      \node (ba)[point] at (2,6) {};    
      \node (bb)[point] at (6,6) {};
      \node (aaa) at (-7, 7) {};
      \node (aab) at (-5, 7) {};
      \node (aba) at (-3, 7) {};
      \node (abb) at (-1, 7) {};
      \node (baa) at (1, 7) {};
      \node (bab) at (3, 7) {};
      \node (bba) at (5, 7) {};
      \node (bbb) at (7, 7) {};
      \node (h21) at (-4, 1) {\huge $\hilb^2(\PP^1)$};
      \node (h11) at (2, -.5) {\huge $\PP^1$};      
      \node (ht+12) at (4, 1.5) {\huge $\hilb^{t+1}(\PP^2)$};
      
      \draw[ultra thick] (o) -- (a);
      \draw[ultra thick] (o) -- (b);
      \draw[ultra thick] (a) -- (aa);
      \draw[ultra thick] (a) -- (ab);
      \draw[ultra thick] (b) -- (ba);
      \draw[ultra thick] (b) -- (bb);
      \draw[ultra thick, dashed] (aa) -- (aaa);
      \draw[ultra thick, dashed] (aa) -- (aab);
      \draw[ultra thick, dashed] (ab) -- (aba);
      \draw[ultra thick, dashed] (ab) -- (abb);
      \draw[ultra thick, dashed] (ba) -- (baa);
      \draw[ultra thick, dashed] (ba) -- (bab);
      \draw[ultra thick, dashed] (bb) -- (bba);
      \draw[ultra thick, dashed] (bb) -- (bbb);
      \draw[thick, -stealth] (h21) .. controls (-5,3) .. (a);
      \draw[thick, -stealth] (h11) .. controls (1,0) .. (o);
      \draw[thick, -stealth] (ht+12) .. controls (5,3) .. (b);

      \node (n2) at (15,-2) {\Huge $\hilbtree_2$};      
      \node (2o)[point] at (15,0) {};
      \node (2a)[point] at (11,4) {};
      \node (2b)[point] at (19,4) {};
      \node (2aa)[point] at (9,6) {};    
      \node (2ab)[point] at (13,6) {};
      \node (2ba)[point] at (17,6) {};    
      \node (2bb)[point] at (21,6) {};
      \node (2aaa) at (8, 7) {};
      \node (2aab) at (10, 7) {};
      \node (2aba) at (12, 7) {};
      \node (2abb) at (14, 7) {};
      \node (2baa) at (16, 7) {};
      \node (2bab) at (18, 7) {};
      \node (2bba) at (20, 7) {};
      \node (2bbb) at (22, 7) {};
      \node (h22) at (11, 1) {\huge $\hilb^2(\PP^2)$};
      \node (h12) at (17, -.5) {\huge $\PP^2$};      
      \node (ht+13) at (19, 1.5) {\huge $\hilb^{t+1}(\PP^3)$};
      \node (g13) at (20, .4) {\huge $ = \mathbb{G}(1,3)$};
      \node (h32) at (11, 5.75) {\Large $\hilb^3(\PP^2)$};
      \node (h2t+13) at (15, 5.75) {\Large $\hilb^{2t+1}(\PP^3)$};
      \node (ht+23) at (19, 5.75) {\Large $\hilb^{t+2}(\PP^3)$};
      \node (g24) at (21, 4.5) {\Large $\mathbb{G}(2,4)$};
      
      \draw[ultra thick] (2o) -- (2a);
      \draw[ultra thick] (2o) -- (2b);
      \draw[ultra thick] (2a) -- (2aa);
      \draw[ultra thick] (2a) -- (2ab);
      \draw[ultra thick] (2b) -- (2ba);
      \draw[ultra thick] (2b) -- (2bb);
      \draw[ultra thick, dashed] (2aa) -- (2aaa);
      \draw[ultra thick, dashed] (2aa) -- (2aab);
      \draw[ultra thick, dashed] (2ab) -- (2aba);
      \draw[ultra thick, dashed] (2ab) -- (2abb);
      \draw[ultra thick, dashed] (2ba) -- (2baa);
      \draw[ultra thick, dashed] (2ba) -- (2bab);
      \draw[ultra thick, dashed] (2bb) -- (2bba);
      \draw[ultra thick, dashed] (2bb) -- (2bbb);
      \draw[thick, -stealth] (h22) .. controls (10,3) .. (2a);
      \draw[thick, -stealth] (h12) .. controls (16,0) .. (2o);
      \draw[thick, -stealth] (ht+13) .. controls (20,3) .. (2b);
      \draw[thick, -stealth] (h32) .. controls (10,6.5) .. (2aa);
      \draw[thick, -stealth] (h2t+13) .. controls (14,6.5) .. (2ab);
      \draw[thick, -stealth] (ht+23) .. controls (18,6.5) .. (2ba);
      \draw[thick, -stealth] (g24) .. controls (22,6.5) .. (2bb);
      
      \node (n3) at (30,-2) {\Huge $\hilbtree_3$};      
      \node (3o)[point] at (30,0) {};
      \node (3a)[point] at (26,4) {};
      \node (3b)[point] at (34,4) {};
      \node (3aa)[point] at (24,6) {};    
      \node (3ab)[point] at (28,6) {};
      \node (3ba)[point] at (32,6) {};    
      \node (3bb)[point] at (36,6) {};
      \node (3aaa) at (23, 7) {};
      \node (3aab) at (25, 7) {};
      \node (3aba) at (27, 7) {};
      \node (3abb) at (29, 7) {};
      \node (3baa) at (31, 7) {};
      \node (3bab) at (33, 7) {};
      \node (3bba) at (35, 7) {};
      \node (3bbb) at (37, 7) {};
      \node (h23) at (26, 1) {\huge $\hilb^2(\PP^3)$};
      \node (h13) at (32, -.5) {\huge $\PP^3$};      
      \node (ht+14) at (34, 1.5) {\huge $\hilb^{t+1}(\PP^4)$};
      
      \draw[ultra thick] (3o) -- (3a);
      \draw[ultra thick] (3o) -- (3b);
      \draw[ultra thick] (3a) -- (3aa);
      \draw[ultra thick] (3a) -- (3ab);
      \draw[ultra thick] (3b) -- (3ba);
      \draw[ultra thick] (3b) -- (3bb);
      \draw[ultra thick, dashed] (3aa) -- (3aaa);
      \draw[ultra thick, dashed] (3aa) -- (3aab);
      \draw[ultra thick, dashed] (3ab) -- (3aba);
      \draw[ultra thick, dashed] (3ab) -- (3abb);
      \draw[ultra thick, dashed] (3ba) -- (3baa);
      \draw[ultra thick, dashed] (3ba) -- (3bab);
      \draw[ultra thick, dashed] (3bb) -- (3bba);
      \draw[ultra thick, dashed] (3bb) -- (3bbb);
      \draw[thick, -stealth] (h23) .. controls (25,3) .. (3a);
      \draw[thick, -stealth] (h13) .. controls (31,0) .. (3o);
      \draw[thick, -stealth] (ht+14) .. controls (35,3) .. (3b);

      \node (n4) at (35,-2) {\Huge $\cdots$};            
    \end{tikzpicture}
  }
  \caption{The forest of Hilbert schemes \label{forest}}
\end{figure}

The first few trees in this arrangement are sketched in
Figure~\ref{forest} with some of the Hilbert schemes corresponding to
vertices indicated.

This framework supplies us with a rough chart of the geography of
Hilbert schemes, highlighting that certain properties of Hilbert
schemes hold in a predictable manner.  For instance, this leads to the
following classification of Hilbert schemes with unique Borel-fixed
points.

\begin{theorem}
  \label{thm:SSSunique}
  Let $\hp$ be an admissible Hilbert polynomial and $c = n - \deg
  \hp$.  The lexicographic point is the unique Borel-fixed point on
  $\hilb^{\hp}(\PP^n)$ if and only if one of the following holds:
  \begin{enumerate}
  \item $b_r > 0$,
  \item $c \ge 2$ and $r \le 2$,
  \item $c = 1$ and $b_1 = b_r$, or
  \item $c = 1$ and $r - s \le 2$, where $b_1 = b_2 = \dotsb = b_{s} >
    b_{s+1} \ge \dotsb \ge b_r$.
  \end{enumerate}
\end{theorem}

\begin{proof}
  Let $I \subset S$ be Borel-fixed with Hilbert polynomial $\hp$
  satisfying one of the conditions.  In
  \cite[Theorem~1.1]{Staal--2020}, the base field $\kk$ is assumed to
  be algebraically closed or have characteristic $0$.  If $\kk$ is
  infinite of characteristic $p > 0$, then $I$ is $p$-Borel (see \S
  \ref{ch:borelfixed}) and base-change to the algebraic closure
  $\overline{\kk} \supseteq \kk$ produces $\overline{I} := I \cdot
  \overline{\kk}[x_0, x_1, \dotsc, x_n]$ with $\hp_{\overline{I}} =
  \hp_I$.  In fact, $\hf_{\overline{I}} = \hf_I$ holds and
  $\overline{I}$ is also $p$-Borel (i.e.\ \cite[Theorem
    1.1.2]{Herzog--Hibi--2011}), so is Borel-fixed and thus
  lexicographic, which implies that $I = \overline{I} \cap S$ is
  lexicographic.
\end{proof}

The \emph{\bfseries lexicographic point} $[L^{\hp}_n] \in
\hilb^{\hp}(\PP^n)$ is the point defined by the saturated lex-segment
ideal $L^{\hp}_n \subseteq S$ with Hilbert polynomial $\hp$.  The
lexicographic point is nonsingular and lies on a unique irreducible
component $\hilb^{\hp}(\PP^n)$ called the \emph{\bfseries
lexicographic} or \emph{\bfseries Reeves--Stillman} component
\cite{Reeves--Stillman--1997}.  It is described in the next section.

\subsection{Lexicographic Ideals}
\label{ch:lexideals}

Lex-segment, or lexicographic, ideals are monomial ideals whose
homogeneous pieces are spanned by maximal monomials in lexicographic
order.  They are Borel-fixed, and in many respects are the most
important monomial ideals.

For any $u = (u_0, u_1, \dotsc, u_n) \in \NN^{n+1}$, let $x^u =
x_0^{u_0} x_1^{u_1} \dotsb x_n^{u_n}$.  The \emph{\bfseries
lexicographic ordering} is the relation $\lexg$ on the monomials in
$S$ defined by $x^u \lexg x^v$ if the first nonzero coordinate of $u -
v \in \ZZ^{n+1}$ is positive, where $u, v \in \NN^{n+1}$.

\begin{example}
  \label{eg:lexordering}
  We have $x_0 \lexg x_1 \lexg \dotsb \lexg x_n$ in lexicographic
  order on $S = \kk[x_0, x_1, \dotsc, x_n]$.  If $n \ge 2$, then $x_0
  x_2^2 \lexg x_1^4 \lexg x_1^3$.
\end{example}

For a homogeneous ideal $I \subseteq S$, lexicographic order gives
rise to two associated monomial ideals.  First, the \emph{\bfseries
lex-segment ideal} for $\hf_I$ is the monomial ideal $L^{\hf_I}_n
\subseteq S$ whose $i$th graded piece is spanned by the $\dim_{\kk}
I_i = \hf_{S}(i) - \hf_I(i)$ largest monomials in $S_i$, for all $i
\in \ZZ$.  The equality $\hf_I = \hf_{L^{\hf_I}_n}$ holds by
definition and $L^{\hf_I}_n$ is a homogeneous ideal of $S$ by \cite[\S
  II]{Macaulay--1927}; see \cite[Proposition
  2.21]{Miller--Sturmfels--2005}.  More importantly, the
\emph{\bfseries (saturated) lexicographic ideal} for $\hp_I$ is the
lex-segment ideal $L^{\hp_I}_n := \bigl( L^{\hf_I}_n :
\mathfrak{m}^{\infty} \bigr)$, where $\mathfrak{m} := \langle x_0,
x_1, \dotsc, x_{n} \rangle \subset S$ is the irrelevant ideal.
Saturation with respect to $\mathfrak{m}$ does not affect the Hilbert
function in large degrees, so $L^{\hp_I}_n$ also has Hilbert
polynomial $\hp_I$.  (The two ideals can coincide for certain $I$,
e.g.\ when $I$ is already saturated and lexicographic, so $I =
L^{\hf_I}_n = L^{\hp_I}_n$.)

Given a finite sequence of nonnegative integers $a_0, a_1, \dotsc,
a_{n-1} \in \NN$, consider the monomial ideal $L(a_0, a_1, \dots,
a_{n-1}) \subset S$ from \cite[Notation~1.2]{Reeves--Stillman--1997}
with generators
\[
\langle x_0^{a_{n-1}+1}, x_0^{a_{n-1}} x_1^{a_{n-2}+1}, \dotsc,
x_{0}^{a_{n-1}} x_{1}^{a_{n-2}} \dotsb x_{n-3}^{a_{2}} x_{n-2}^{a_{1}
  + 1}, x_{0}^{a_{n-1}} x_{1}^{a_{n-2}} \dotsb x_{n-2}^{a_{1}}
x_{n-1}^{a_{0}} \rangle .
\]
These are the saturated lexicographic ideals.

\begin{lemma}
  \label{lem:lexmingens}
  Let $\hp(t) = \sum_{j=1}^{r} \binom{t + b_j - j + 1}{b_j}$ have
  Gotzmann partition $(b_1, b_2, \dotsb, b_r)$, let $a_j$ be the
  number of parts equal to $j$, and let $n \in \NN$ satisfy $n > d :=
  \deg \hp$.
  \begin{enumerate}
  \item We have
    \begin{align*}
      L^{\hp}_n &= L(a_0, a_1, \dotsc, a_{n-1}) \\
                &= \langle x_0, x_1, \dotsc , x_{n-d-2},
                  x_{n-d-1}^{a_{d} + 1}, \\
                &\relphantom{= \langle} x_{n-d-1}^{a_{d}}
                  x_{n-d}^{a_{d-1} + 1}, \dotsc ,
                  x_{n-d-1}^{a_{d}} x_{n-d}^{a_{d-1}} \dotsb
                  x_{n-3}^{a_{2}} x_{n-2}^{a_{1} + 1}, 
                  x_{n-d-1}^{a_{d}}
                  x_{n-d}^{a_{d-1}} \dotsb x_{n-2}^{a_{1}}
                  x_{n-1}^{a_{0}} \rangle .
    \end{align*}

  \item If there is an integer $0 \le \ell \le d-1$ such that $a_j =
    0$ for all $j \le \ell$, and $a_{\ell + 1} > 0$, then the minimal
    monomial generators of $L^{\hp}_n$ are given by $m_1, m_2, \dotsc,
    m_{n - \ell - 1}$, where
    \begin{align*}
      m_i &= x_{i-1}, \text{ for all } 1 \le i \le n-d-1, \\
      m_{n-d + k} &= \left( \prod_{j=0}^{k-1} x_{n-d-1 + j}^{a_{d -
          j}} \right) x_{n-d-1 + k}^{a_{d-k} + 1}\,, \text{ for all } 0 \le k \le d -
        \ell -2, \text{ and } \\
      m_{n-\ell-1} &= \prod_{j=0}^{d-\ell-1} x_{n-d-1 + j}^{a_{d - j}}.
    \end{align*}
    If $a_0 \ne 0$, then the minimal monomial generators are those
    listed in (i).
  \end{enumerate}
\end{lemma}

\begin{proof} 
  See \cite[Theorem~2.3]{Moore--Nagel--2014} and
  \cite[Lemma~2.10]{Staal--2020}.
\end{proof}

Importantly, Lemma~\ref{lem:lexmingens} shows that all sequences
$(a_0, a_1, \dotsc, a_{n-1})$ of nonnegative integers determine a
lexicographic ideal.

The two operations $\plus$ and $\lift$ on admissible Hilbert
polynomials have analogues on lexicographic ideals.  For any ideal $I
\subseteq S$, we use $\lift(I) := I \cdot S[x_{n+1}]$ to denote the
lifted ideal.

\begin{corollary} 
  \label{cor:lexoperations}
  Let $\hp$ be admissible, $n > \deg \hp$ an integer, and $L^{\hp}_n =
  L(a_0, a_1, \dotsc, a_{n-1})$.  We have $\plus(L^{\hp}_n) :=
  L^{\plus(\hp)}_n = L(a_0 + 1, a_1, a_2, \dotsc, a_{n-1})$ and $\lift
  \bigl( L^{\hp}_n \bigr) = L^{\lift(\hp)}_{n+1} = L(0, a_0, a_1,
  \dotsc, a_{n-1})$.
\end{corollary}

\begin{proof}
  These both follow directly from Lemma~\ref{lem:lexmingens}(i).
  Cf.\ Proposition/Definition~\ref{propdef:hilbtree}.
\end{proof}

\section{Borel-fixed Ideals}
\label{ch:borelfixed}

Here we review some essential properties of Borel-fixed ideals.  An
ideal $I \subseteq S$ is \emph{\bfseries Borel-fixed} if it is fixed
by the action $\gamma \cdot x_j = \sum_{i=0}^{n} \gamma_{ij} x_i$ of
upper triangular matrices $\gamma \in \operatorname{GL}_{n+1}(\kk)$.
Pardue \cite[II]{Pardue--1994} gives the following combinatorial
criterion for $I$ to be Borel-fixed when $\kk$ is infinite of
characteristic $p > 0$: $I$ is monomial and for each monomial $x^u \in
I$, if $x_i \lexg x_j$, then $x_j^{-k} x_i^k x^u \in I$ holds, for all
$k \le_p u_j$.  Here, $k \le_p \ell$ means that in the base-$p$
expansions $k = \sum_i k_i p^i, \ell = \sum_i \ell_i p^i$, we have
$k_i \le \ell_i$, for all $i$.  This criterion makes sense over any
field and an ideal satisfying it is called \emph{\bfseries $p$-Borel}.

\begin{example}
  The ideal $I := \llrr{x_0^2 + x_0x_1 + x_1^2} \in \FF_2[x_0, x_1,
    x_2]$ is Borel-fixed with Hilbert polynomial $\hp_I = 2t+1$.  This
  shows that Theorem \ref{thm:SSSunique} does not hold over finite
  fields, as the Gotzmann partition of $\hp_I$ is $(1,1)$.
\end{example}

When $\kk$ has characteristic $0$, Pardue's criterion reduces to the
following property (which again makes sense over any field).  A
monomial ideal $I \subseteq S$ is \emph{\bfseries strongly stable}, or
\emph{\bfseries $0$-Borel}, if, for all monomials $m \in I$, all $x_j
| m$, and all $x_i \lexg x_j$, we have $x_j^{-1} x_i m \in I$.  A
strongly stable ideal is always Borel-fixed, but the converse is not
true in positive characteristic.  A Borel-fixed ideal that is not
strongly stable is called a \emph{\bfseries nonstandard} Borel-fixed
ideal.

For a monomial $m$, let $\max m$ be the maximum index $j$ such that
$x_j | m$, and $\min m$ be the minimum such index.

\begin{example}
  The monomial ideal $I = \langle x_0^2, x_0 x_1, x_1^2 \rangle
  \subset S$ is strongly stable when $n \ge 1$.  The ideal $I :=
  \langle x_0^p, x_1^p \rangle \subset S$ is a nonstandard Borel-fixed
  ideal if $\kk$ has characteristic $p > 0$; it is $p$-Borel in any
  characteristic.  The monomial $m = x_1^5 x_2 x_7^2 \in \kk[x_0, x_1,
    \dotsc, x_{13}]$ satisfies $\max m = 7$ and $\min m = 1$.
\end{example}

For a monomial ideal $I$, let $G(I)$ denote its minimal set of
monomial generators.  Whether $I$ is Borel-fixed can be checked on
monomials in $G(I)$; see \cite[\S 15.9.3]{Eisenbud--1995}.

Saturated strongly stable ideals are generated by expansions and
lifting.  Let $I \subseteq S$ be a saturated strongly stable ideal.  A
generator $g \in G(I)$ is \emph{\bfseries expandable} if there are no
elements of $G(I)$ in the set $\left\{ x_i^{-1} x_{i+1} g \mid x_i
\text{ divides } g \text{ and } 0 \le i < n-1 \right\}$.  The
\emph{\bfseries expansion} of $I$ at an expandable generator $g$ is
the monomial ideal $I' \subseteq S$ with minimal generators
\[ 
G(I') := (G(I) \setminus \{ g \}) \cup \{ g x_{j} \mid \max g \le j
\le n-1 \};
\]
see \cite[Definition~3.4]{Moore--Nagel--2014}.  The monomial $1 \in
\langle 1 \rangle$ is vacuously expandable with expansion
$\mathfrak{m} \subset S$.  The expansion of a saturated strongly
stable ideal is again strongly stable, by definition.  It is also
saturated, as we see in Lemma~\ref{lem:HpolyBorel}.

\begin{example}
  \label{eg:lexexpansion}
  Let $L^{\hp}_n = L(a_0, a_1, \dotsc, a_{n-1}) = \langle m_1, m_2,
  \dotsc, m_{n - \ell - 1} \rangle$, as in Lemma~\ref{lem:lexmingens}.
  The expansion at $m_{n-\ell-1}$ is easily verified to be $L(a_0 + 1,
  a_1, a_2, \dotsc, a_{n-1}) = \plus \bigl( L^{\hp}_n \bigr)$.
\end{example}

For a Borel-fixed ideal $I \subseteq S$, let $\nabla(I) \subseteq
\kk[x_0, x_1, \dotsc, x_{n-1}]$ be the image of $I$ under the mapping
$S = \kk[x_0, x_1, \dotsc, x_{n}] \to \kk[x_0, x_1, \dotsc, x_{n-1}]$,
defined by $x_j \mapsto x_j$ for $0 \le j \le n-2$ and $x_j \mapsto 1$
for $n-1 \le j \le n$.  The following lemma shows how this relates to
saturation and generalizes Corollary~\ref{cor:lexoperations} to
Borel-fixed ideals.

\begin{lemma}
  \label{lem:HpolyBorel}
  Let $I \subseteq S$ be a Borel-fixed ideal over an infinite field.
  \begin{enumerate}
  \item We have $(I : \mathfrak{m}_k^{\infty}) = (I : x_k^{\infty})$,
    where $\mathfrak{m}_k := \langle x_0, x_1, \dotsc, x_{k} \rangle
    \subset S$, for all $0 \le k \le n$.
    
  \item If $I$ is saturated, then $\hp_{\nabla(I)} = \nabla (\hp_I)$
    holds, where $\nabla(\hp)(t) := \hp(t) - \hp(t-1)$.

  \item If $I$ is saturated, then there exists $j \in \NN$ such that
    $\hp_{\lift(I)} = \plus^j \lift (\hp_I)$.

  \item For a saturated strongly stable $I$ with expansion $I'$, we
    have $\hp_{I'} = \plus (\hp_I) := 1 + \hp_I$.
  \end{enumerate}
\end{lemma}

\begin{proof} $\;$
  \begin{enumerate}
  \item See \cite[II, Proposition~9]{Pardue--1994}.

  \item Let $J = I \cap \kk[x_0, x_1, \dotsc, x_{n-1}]$.  Because $I$
    is saturated, $x_n$ is not a zero-divisor and
  \[ 
  0 \longrightarrow \left( S/I \right)(-1) \longrightarrow S/I
  \longrightarrow \kk[x_0, x_1, \dotsc, x_{n-1}] / J \longrightarrow 0
  \]
  is a short exact sequence.  The Hilbert function of $J$ now
  satisfies $\hf_{J}(i) = \hf_I(i) - \hf_I(i-1)$, for all $i \in \ZZ$.
  Saturating $J$ with respect to $\langle x_0, x_1, \dotsc, x_{n-1}
  \rangle$ gives $\nabla(I)$, so we find that $\hp_{\nabla(I)}(t) =
  \hp_I(t) - \hp_I(t-1) = \nabla(\hp_I)(t)$.

  \item The ideal $\lift(I)$ is saturated by (i) and Borel-fixed by
    Pardue's criterion, with no elements of $G(\lift(I))$ divisible by
    $x_n$.  Thus, $\nabla \bigl( \lift ( I ) \bigr) = I$.  Then
    $\nabla (\hp_{\lift ( I )}) = \hp_I$ follows by (ii), so $\lift
    \nabla (\hp_{\lift ( I )}) = \lift (\hp_I)$.  Now
    \cite[Lemma~2.5]{Staal--2020} shows $\hp_{\lift(I)} - \lift \nabla
    (\hp_{\lift(I)}) \in \NN$.
    
  \item See \cite[Lemma~3.15]{Moore--Nagel--2014}. \qedhere
  \end{enumerate}
\end{proof}

Lemma~\ref{lem:HpolyBorel} highlights important properties that are
used together with expansions to generate all saturated strongly
stable ideals.  However, we also wish to understand nonstandard
Borel-fixed ideals.  Fortunately, other than (iv), these properties
hold for all Borel-fixed ideals and are useful in
Section~\ref{ch:classify}.

An algorithm for generating saturated strongly stable ideals was first
developed by Reeves in \cite{Reeves--1992}; we follow the approach of
\cite{Moore--Nagel--2014}.  The heart of Reeves' algorithm is the
following. 

\begin{theorem}
  \label{thm:algorithm}
  If $I$ is a saturated strongly stable ideal of codimension $c$, then
  there is a finite sequence $I_{(0)}, I_{(1)}, \dotsc, I_{(i)}$ such
  that $I_{(0)} = \langle 1 \rangle = \kk[x_0, x_1, \dotsc, x_c]$,
  $I_{(i)} = I$, and $I_{(j)}$ is an expansion or lifting of
  $I_{(j-1)}$, for all $1 \le j \le i$.
\end{theorem}

\begin{proof}
  See \cite[Theorem~4.4]{Moore--Nagel--2014}.
\end{proof}

When $I \neq \langle 1 \rangle$, the first step is always to expand
$\langle 1 \rangle$ to $\langle x_0 , x_1, \dotsc, x_{c-1} \rangle$.
We suppress this step, as we assume $\hp \neq 0$.  The sequences of
expansions and liftings are not generally unique, but
Corollary~\ref{cor:lexoperations} shows that they are for
lexicographic ideals.  Theorem~\ref{thm:algorithm} leads to the
following algorithm; see \cite[Appendix~A]{Reeves--1992} and \cite[\S
  4]{Moore--Nagel--2014}.

\begin{algorithm}
  \label{alg:SSS}
  \qquad
  \begin{tabbing}
    \qquad \= \qquad \= \qquad \= \kill
    \textrm{Input:} an admissible Hilbert polynomial $\hp \in \QQ[t]$ and
    $n \in \NN$ satisfying $n > \deg \hp$ \\

  \textrm{Output:} all saturated strongly stable ideals with Hilbert
  polynomial $\hp$ in $\kk[x_0, x_1, \dotsc, x_n]$ \\[0.4em]
    
    $j = 0$; $d = \deg \hp$; \\[0.1em] 

    $\hq_0 = \nabla^d(\hp)$; $\hq_{1} = \nabla^{d-1} (\hp)$; \ldots
    ; $\hq_{d-1} = \nabla^{1} (\hp)$; $\hq_d = \hp$; \\[0.1em]

    $\mathcal{S} = \{ \langle x_0, x_1, \dotsc, x_{n-d-1} \rangle \}$,
  where $\langle x_0, x_1, \dotsc, x_{n-d-1} \rangle \subset\kk[x_0,
    x_1, \dotsc, x_{n-d}]$; \\[0.2em]

    \textrm{WHILE} $j \le d$ \textrm{DO} \\

    \>  $\mathcal{T} = \varnothing$; \\

    \>  \textrm{FOR} $J \in \mathcal{S}$, considered as an ideal in
    $\kk[x_0, x_1, \dotsc, x_{n-d+j}]$ \textrm{DO} \\

    \> \>    \textrm{IF} $\hq_j - \hp_J \ge 0$ \textrm{THEN} \\

    \> \> \> compute all sequences of $\hq_j - \hp_J$ expansions that
    begin with $J$; \\
    \> \> \> $\mathcal{T} = \mathcal{T} \cup \text{ the resulting set of
      sat.\ str.\ st.\ ideals with Hilbert polynomial } \hq_j$; \\

    \>  $\mathcal{S} = \mathcal{T}$; \\
    \> $j = j+1$; \\

    \textrm{RETURN} $\mathcal{S}$
  \end{tabbing}
\end{algorithm}

\begin{proof}
  See \cite[Algorithm~4.6]{Moore--Nagel--2014}.  Here, $\mathcal{S}$
  is reset at the $j$th step to Moore's $\mathcal{S}^{(d-j)}$.
\end{proof}

\section{Classifying Hilbert Schemes with Two Borel-fixed Points}
\label{ch:classify}

We now prove that the classification of Hilbert schemes with exactly
two Borel-fixed points given in \cite{Ramkumar--2019} holds in all
characteristics, after one minor adjustment.  There are two main tools
used to derive the classification in characteristic $0$, namely,
Reeves' algorithm and Theorem~\ref{thm:SSSunique}.  Reeves' algorithm
does not produce nonstandard Borel-fixed ideals, so we must alter our
approach.  Fortunately, the recursive properties of
Lemma~\ref{lem:HpolyBorel} together with the classification in
Theorem~\ref{thm:SSSunique} are enough to show that all Borel-fixed
points are strongly stable on the Hilbert schemes we consider.

We denote $R := \kk[x_0, x_1, \dotsc, x_{n-1}] \supset \mathfrak{n}_k
:= \langle x_0, x_1, \dotsc, x_k \rangle$, for $0 \le k \le n-1$.  As
before, we also denote $d := \deg \hp$ and $n := d + c$.  We start
with some basic facts and only study the case $c \ge 2$, by the
reduction in \cite[\S 2]{Ramkumar--2019}.

\begin{lemma}
  \label{lem:twisted}
  The Hilbert scheme $\hilb^{3t+1}(\PP^n)$ has $3$ Borel-fixed points,
  for all $n \ge 3$.
\end{lemma}

\begin{proof}
  Suppose $I \subset S$ is a saturated Borel-fixed ideal with Hilbert
  polynomial $\hp = 3t+1$.  So $I$ is generated in degrees up to the
  Gotzmann number $r = 4$.  The ideal $\nabla(I)$ is saturated and
  Borel-fixed with Hilbert polynomial $\hp_{\nabla(I)} = 3$, by
  Lemma~\ref{lem:HpolyBorel}.  Thus, $\nabla(I)$ equals either
  $\mathfrak{n}_{c-3} + \langle x_{c-2}^2, x_{c-2} x_{c-1}, x_{c-1}^2
  \rangle$ or $\mathfrak{n}_{c-2} + \langle x_{c-1}^3 \rangle$, by
  Theorem~\ref{thm:main2}(i)(a) (here $n-1 = c$).
  Lemma~\ref{lem:HpolyBorel}(i) shows that $\nabla(I) = (J :
  x_{c}^{\infty})$, where $J := I \cap R$.  Examining generators, this
  means $\lift(\nabla(I)) = (I : x_{n-1}^{\infty}) \supseteq I$.  When
  $\nabla(I) = \mathfrak{n}_{c-3} + \langle x_{c-2}^2, x_{c-2}
  x_{c-1}, x_{c-1}^2 \rangle$ holds, the lift $\lift(\nabla(I))$ has
  Hilbert polynomial $3t+1$ and must equal $I$.  When $\nabla(I)$ is
  lexicographic, we get $I \subset L := \mathfrak{m}_{n-3} + \langle
  x_{n-2}^3 \rangle \subset S$.  In particular, $I_4 \subset L_4$ has
  codimension $1$.  The monomials $x_0 x_n^3, x_1 x_n^3, \dotsc,
  x_{n-3} x_n^3, x_{n-2}^3 x_n$ cannot all belong to $I_4$, as $I$ is
  saturated, so one must be excluded.  By Pardue's criterion, the
  possibilities are to exclude $x_{n-3} x_n^3$ or $x_{n-2}^3 x_n$.  If
  we exclude the first, then $I$ contains, and thus equals, the
  expansion of $L$ at $x_{n-3}$.  If we exclude the second, then $I$
  equals $L^{\hp}_n$.
\end{proof}

The proof of Lemma~\ref{lem:twisted} demonstrates key aspects of the
strategy for proving Theorem~\ref{thm:main1}.  That is, there are
limited possibilities for what $\nabla(I)$ can be, which restricts
what $(I : x_{n-1}^{\infty})$ can be, given a saturated Borel-fixed
ideal $I$.  Then, because we are seeking Hilbert schemes with few
Borel-fixed points, the Hilbert polynomials of $I$ and $(I :
x_{n-1}^{\infty})$ can only differ by a small constant.  This, in
turn, means that $I$ can be effectively derived from the ideal $(I :
x_{n-1}^{\infty})$ using Pardue's criterion.

We need one last lemma before proving the main theorem.

\begin{lemma}
  \label{lem:twoborels}
  Let $n = d + c$, for $c \ge 2$.  If $\hilb^{\hp}(\PP^n)$ has exactly
  two Borel-fixed points, then $L^{\hp -1}_n$ is either
  $\mathfrak{m}_{c-2} + \langle x_{c-1}^a \rangle$, for some $a >1$,
  or $\mathfrak{m}_{c-2} + x_{c-1}^a \langle x_{c-1}, x_c, \dotsc,
  x_{n'-1} \rangle$, for some $a > 0$ and $c < n' \le n$.
\end{lemma}

\begin{proof}
  By Theorem~\ref{thm:SSSunique}, $\hp -1$ is admissible.  Then
  Algorithm~\ref{alg:SSS} implies that $L^{\hp -1}_n$ is generated in
  at most two degrees; see \cite[Proposition~4.4]{Staal--2020}.  If
  $L^{\hp -1}_n$ is generated in a single degree, then it must equal
  $\mathfrak{m}_{c-1} \subset S$, by Lemma~\ref{lem:lexmingens}.  But
  Theorem~\ref{thm:SSSunique} then implies $\hilb^{\hp}(\PP^n)$ has a
  unique Borel-fixed point.  Lemma~\ref{lem:lexmingens}(ii) now leaves
  the provided options.
\end{proof}

Using this preliminary classification, we can now prove the main
theorem.

\begin{theorem}
  \label{thm:main2}
  Let $n = d + c$, where $d = \deg \hp$ and $c \ge 2$.  The Hilbert
  scheme $\hilb^{\hp}(\PP^n)$ has exactly two Borel-fixed points if
  and only if one of the following conditions holds:
  \begin{enumerate}[]
  \item
    \begin{enumerate}[label={(\alph*)}]
    \item $b_1 = 0$ and $r = 3$,

    \item[(a')] $b_1 = 0$ and $r = 4$, if $n = 2$ and $\kk$ does not
      have characteristic $2$,
    
    \item $b_1 = \dotsb = b_{r-1} = 1$ and $b_r = 0$, for $r-1 \neq 1,
      3$,

    \item $b_1 = \dotsb = b_{r-1} \ge 2$ and $b_r = 0$, for $r-1 \neq
      1$,
    \end{enumerate}

  \item
    \begin{enumerate}[label={(\alph*)}]
    \item $b_1 > b_2 = 0$ and $r = 3$,
      
    \item $b_1 = \dotsb = b_{r-2} > b_{r-1} = 1$ and $b_r = 0$, for
      $r-2 \neq 2$, and

    \item $b_1 = \dotsb = b_{r-2} > b_{r-1} \ge 2$ and $b_r = 0$.
    \end{enumerate}
  \end{enumerate}
\end{theorem}

\begin{proof}
  Let $I \subset S$ be saturated, Borel-fixed, and have Hilbert
  polynomial $\hp$ with Gotzmann partition $b := (b_1, b_2, \dotsc,
  b_r)$.  If the corresponding point $[I]$ lies on a
  Hilbert scheme having exactly two Borel-fixed ideals, then one of
  the following holds:
  \begin{enumerate}
  \item $a = r-1 > 1$ and $b_1 = \dotsb = b_{a} \ge b_r = 0$, or
  \item $a = r-2 > 0$ and $b_1 = \dotsb = b_{a} > b_{a+1} \ge b_r =
    0$,
  \end{enumerate}
  by Lemma~\ref{lem:twoborels}.  When $b$ does not satisfy (i) or
  (ii), either $I$ is lexicographic, by Theorem~\ref{thm:SSSunique},
  or there are at least three saturated strongly stable ideals with
  Hilbert polynomial $\hp$, by Reeves' algorithm, and these are
  Borel-fixed.  Thus, we suppose $I$ satisfies either (i) or (ii).
  The analyses of several cases are analogous to the proof of
  Lemma~\ref{lem:twisted}.

  \vspace{0.25em}
  
  \paragraph{(i)(a)}
  Here $b = (0, 0, 0)$, i.e.\ $\hp = 3$ and $n = c$.  We know that $I$
  is generated in degrees up to $r = 3$.  Suppose $n = 2$, so $I
  \subset S = \kk[x_0, x_1, x_2]$.  Then $x_1^N \in G(I)$, for some
  minimal $N \le 3$, by \cite[II, Corollary~8]{Pardue--1994}.  If $N =
  1$, then $x_0 \in I$ too, which implies $\hp = 1$.  If $N = 3$, then
  $x_1^2 x_2^{j-2}, x_1 x_2^{j-1}, x_2^j \notin I_j$, for all $j \gg
  0$, so $x_0 \in I$ and $I = \langle x_0, x_1^3 \rangle$ is
  lexicographic.  If $N = 2$, then $x_0^2 \in I$.  Thus one of the
  monomials $x_0 x_1 x_2, x_0 x_2^2, x_1 x_2^2, x_2^3$ must belong to
  $I_3$.  As $I$ is saturated, we can only have $x_0 x_1 x_2 \in I$
  and so $I = \langle x_0^2, x_0 x_1, x_1^2 \rangle$ is strongly
  stable.  Now suppose $n > 2$ and $x_0 \notin I$.  Then Pardue's
  criterion implies $x_0 x_n^{j-1}, x_1 x_n^{j-1}, \dotsc, x_n^j
  \notin I_j$, for all $j \gg 0$, giving a contradiction.  So $x_0 \in
  I$ and we are finished, by induction on $n$.
  
  \vspace{0.25em}
  
  \paragraph{(i)(a')}
  Here $b = (0,0,0,0)$, so $\hp = 4$ and $n = c = 2$.  (When $n > 2$,
  there are three saturated strongly stable ideals with Hilbert
  polynomial $\hp = 4$.)  As in (i)(a), we can see that $G(I)$
  contains one of $x_1^4, x_1^3$, or $x_1^2$.  We summarize: if $x_1^4
  \in G(I)$, then $I = \langle x_0, x_1^4 \rangle$ must be
  lexicographic; if $x_1^3 \in G(I)$, then we must have $I = \langle
  x_0^2, x_0 x_1, x_1^3 \rangle$; and if $x_1^2 \in G(I)$, then we can
  only have $I = \langle x_0^2, x_1^2 \rangle$, which is Borel-fixed
  if and only if $\kk$ has characteristic $2$.  If $\hp = r > 4$, then
  Algorithm~\ref{alg:SSS} generates at least $3$ saturated strongly
  stable ideals.

  \vspace{0.25em}
  
  \paragraph{(i)(b)}
  Here $b = (1, 1, \dotsc, 1, 0)$ with $a = r-1 \neq 3$, so that
  $\hp(t) = at + 2 - \frac{(a-1)(a-2)}{2}$ and $\nabla(\hp) = a$.  We
  study $I$ by lifting $\nabla(I)$.  If $a = 2$, then the result
  follows as in Lemma~\ref{lem:twisted}, so let $a \ge 4$.  The goal
  is to show that if $J \subset R$ is a nonlexicographic saturated
  Borel-fixed ideal with Hilbert polynomial $a$, then $\hp_{\lift(J)}
  - \hp > 0$.  This holds for strongly stable $J$, by
  \cite[Theorem~2.15]{Ramkumar--2019}, so we assume that $J$ is
  nonstandard.  As $J$ is nonlexicographic, we have $x_{c-1}^N \in
  G(J)$, for some $2 \le N < a$, and also $x_{c-2} \notin J$.
  Consider the minimal strongly stable ideal $J^0 \supset J$; we
  obtain $J^0$ from $J$ by
  including all ``missing'' monomials of the form $(x_{j_1} x_{j_2}
  \cdots x_{j_s})^{-1} x_{i_1} x_{i_2} \cdots x_{i_s} h \notin J$,
  where $h \in G(J)$, $i_k < j_k$, and $x_{j_k} | h$, for $1 \le k \le
  s$.  Let $a' := \hp_{J^0}$, so that $a' < a$.  Because $x_{c-2}
  \notin J$, it follows that $J^0$ is not lexicographic, which further
  shows $N < a'$.

  Fix $j \gg 0$.  We have $x_{c-1}^{N-1} x_c^{j-N+1}, x_{c-1}^{N-2}
  x_c^{j-N+2}, \dotsc, x_c^{j} \notin J^0_j$.  As $J^0$ is strongly
  stable and contains $x_{c-1}^N$, $R_j / J^0_j$ has basis
  $(x_{c-1}^{N-i} x_c^{j-N+i}, m_{k} x_c^{j - \deg m_k} \mid 1 \le i
  \le N, 1 \le k \le a'-N)$, where the monomials $m_1, m_2, \dotsc,
  m_{a'-N}$ are not divisible by $x_c$ and $\deg m_k < N$, for all $1
  \le k \le a'-N$.  To describe $R_j / J_j$, we need a further $a-a'$
  monomials from $J^0_j \setminus J_j$.  First, as $J^0$ is strongly
  stable, it is generated in degrees up to $N$.  The regularity of $J$
  is at most $a$, so there are $a-a'$ monomials in $J^0_a \setminus
  J_a$.  Each has a unique factorization $h = g h'$, where $g \in
  G(J^0) \setminus J$ and $\max g \le \min h'$; see
  \cite[Lemma~2.11]{Miller--Sturmfels--2005}.  We need to alter these
  monomials to show their effect on the Hilbert function
  $\hf_{\lift(J)}$.  Fix $g \in G(J^0) \setminus J$ for the moment and
  let $h_0, h_1, \dotsc, h_{i_g}$ be the monomials of $J^0_a \setminus
  J_a$ whose factorization involves $g$.  Write $h_i := g p_i x_c^{a -
    \deg g - \deg p_i}$, where $x_c$ does not divide $p_i$.  If $i_g >
  0$, order these such that $h_0 := g x_c^{a - \deg g}$ and $0 < \deg
  p_1 \le \deg p_2 \le \dotsb \le \deg p_{i_g}$.  Let $p_1'$ be any
  monomial dividing $p_1$ with $\deg p_1' = 1$.  Then $f_1 := g p_1'
  \in J^0 \setminus J$ holds.  If $i_g > 1$, then do the following for
  each $1 < i \le i_g$: if $\deg p_i \le i$, then set $p_i' := p_i$,
  otherwise let $p_i'$ divide $p_i$ with $\deg p_i' = i$.  Then we
  have $f_i := g p_i' \in J^0 \setminus J$.  Moreover, $f_i$ differs
  from $f_0 := g, f_1, \dotsc, f_{i-1}$.  To see this, assume $f_i =
  f_k$, for $k < i$, so that $p_i' = p_k'$.  This implies $\deg p_i' =
  \deg p_k' \le k < i$, so $p_i' = p_i$.  But then $p_i$ divides $p_k$
  so that $p_i = p_k$, a contradiction.  We now write $f^g_i$ to
  distinguish the monomials obtained from a fixed $g$.  By uniqueness
  of the factorizations, we have $f^g_i \neq f^{g'}_{i'}$ when $g \neq
  g'$.

  Lifting,
  we see that $S_j \setminus \lift(J)_j$ contains the monomials
  $x_{n-2}^k \langle x_{n-1}, x_n \rangle^{j-k}$, for $0 \le k \le
  N-1$, along with $m_k \langle x_{n-1}, x_n \rangle^{j-\deg m_k}$,
  for $1 \le k \le a'-N$, and $f^g_k \langle x_{n-1}, x_n
  \rangle^{j-\deg f^g_k}$, for $g \in G(J^0) \setminus J$ and $0 \le k
  \le i_g$.  This implies
  \begin{align*}
    \hf_{\lift(J)}(j) &\ge \sum_{k = 0}^{N-1} (j-k+1) + \sum_{k =
      1}^{a'-N} (j - \deg m_k +1) + \sum_{g \in G(J^0) \setminus J}
    \sum_{k=0}^{i_g} (j - \deg f^g_k +1) \\ &\ge \sum_{i = 0}^{N-1} (j
    - i + 1) + \sum_{i = N}^{a'-1} (j - N + 1 +1) + \sum_{i=a'}^{a-1}
    (j - i + 1 + 1) \\ &\ge aj +1 - \frac{(a-1)(a-2)}{2} + 2 \ge
    \hp(j) + 1,
  \end{align*}
  which shows that $\hp_{\lift(J)} - \hp > 0$, as desired.

  \vspace{0.25em}
  
  \paragraph{(i)(c)}
  Here $b = (d, d, \dotsc, d, 0)$, $d \ge 2$, and $a = r-1 \ge 2$.
  Lemma~\ref{lem:HpolyBorel} and Theorem~\ref{thm:SSSunique} imply
  $\nabla(I) = L^{\nabla(\hp)}_{n-1} = \mathfrak{n}_{c-2} + \langle
  x_{c-1}^a \rangle \subset R$.  It follows that $L := L^{\hp -1}_n =
  (I : x_{n-1}^{\infty})$ and $I_{a+1} \subset L_{a+1}$ has
  codimension $1$.  The monomials $x_0 x_n^a, x_1 x_n^a, \dotsc,
  x_{c-2} x_n^a, x_{c-1}^a x_n \in L_{a+1}$ cannot all belong to $I$.
  Pardue's criterion implies either $x_{c-2} x_n^a \notin I_{a+1}$ or
  $x_{c-1}^a x_n \notin I_{a+1}$.  In the first case, $I$ is the
  expansion of $L$ at $x_{c-2}$ and in the second case, $I$ is the
  expansion of $L$ at $x_{c-1}^a$.
  
  \vspace{0.25em}

  \paragraph{(ii)(a)}
  Here $b = (d, 0, 0)$ with $d > 0$, so $\hp = \binom{t + d}{d} + 2$.
  Lemma~\ref{lem:HpolyBorel} and Theorem~\ref{thm:SSSunique} imply
  that $\nabla(I) = L^{\nabla(\hp)}_{n-1} = \mathfrak{n}_{c-1}$.  This
  implies $L := L^{\hp - 2}_n = \mathfrak{m}_{c-1} = (I :
  x_{n-1}^{\infty}) \supset I$, so $I_3$ has codimension $2$ in $L_3$.
  As $I$ is saturated, if $x_{c-1} x_n^2 \in I$, then $I \supseteq L$
  holds.  So $x_{c-1} x_n^2 \notin I$.  Now suppose $I$ contains
  $x_{c-2} x_n^2$.  Then $\mathfrak{m}_{c-2} \subset I$ holds.  If we
  also have $x_{c-1} x_{n-1} x_n \in I$, then $I \supseteq L^{\hp
    -1}_n$ must hold, which is nonsense.  Thus, if $x_{c-2} x_n^2 \in
  I$, then $I_3$ contains all monomials in $L_3$ except $x_{c-1}
  x_{n-1} x_n, x_{c-1} x_n^2$.  This implies that $I = L^{\hp}_n$.
  Now suppose $I$ does not contain $x_{c-2} x_n^2$.  Then $I$ equals
  the expansion of $L^{\hp -1}_n$ at $x_{c-2}$.  If $a \neq 1$ and
  $b_{a+1} = 0$, then Reeves' algorithm generates a third saturated
  Borel-fixed ideal; see \cite[Theorem~2.16]{Ramkumar--2019}.

  \vspace{0.25em}

  \paragraph{(ii)(b)}
  Here $b = (d, d, \dotsc, d, 1, 0)$ with $d > 1$ and $a = r-2 \ne 2$.
  If $a = 1$, then $\nabla(I) = L^{\nabla(\hp)}_{n-1}$, by
  Lemma~\ref{lem:HpolyBorel} and Theorem~\ref{thm:SSSunique}.  Setting
  $L := \mathfrak{m}_{c-2} + x_{c-1} \langle x_{c-1}, x_c, \dotsc,
  x_{n-2} \rangle$, we see that $I_3 \subset L_3$ is a codimension $1$
  subspace.  The elements
  \[
  x_0 x_n^2, x_1 x_n^2, \dotsc, x_{c-2} x_n^2, x_{c-1}^2 x_n, x_{c-1}
  x_c x_n, \dotsc, x_{c-1} x_{n-2} x_{n}
  \]
  cannot all belong to $I_3$.  By Pardue's criterion, either $x_{c-2}
  x_n^2$ or $x_{c-1} x_{n-2} x_{n}$ does not belong.  If $x_{c-2}
  x_n^2 \notin I$, then $I$ is the expansion of $L$ at $x_{c-2}$; if
  $x_{c-1} x_{n-2} x_{n} \notin I$, then $I$ is the expansion of $L$
  at $x_{c-1} x_{n-2}$.

  Suppose $a = 3$.  If $d = 2$, then there are three saturated
  Borel-fixed ideals with Hilbert polynomial $\nabla(\hp) = 3t+1$ and
  codimension $c$, by Lemma~\ref{lem:twisted}.  These are all strongly
  stable and one can directly verify (by reducing to $c=2$) that the
  lifts of the nonlexicographic ones have Hilbert polynomial $\hp +1$.
  Thus, $\nabla(I) = L^{3t+1}_{n-1}$ holds and $I_5 \subset L_5$ has
  codimension $1$, where $L = \mathfrak{m}_{c-2} + \langle x_{c-1}^4,
  x_{c-1}^3 x_c \rangle$.  Again, Pardue's criterion tells us that $I$
  is an expansion of $L$.  If $d \ge 3$ or $a \ge 4$, then $\nabla(I)$
  falls into cases (i)(b)-(c) and must equal either
  \[
  \mathfrak{n}_{c-2} + x_{c-1}^a \llrr{x_{c-1}, x_c, \ldots, x_{n-2}}
  \quad \text{ or } \quad \mathfrak{n}_{c-3} + x_{c-2} \llrr{x_{c-2},
    x_{c-1}, \ldots, x_{n-2}} + \llrr{x_{c-1}^a}
  \]
  in $R$.  Letting $L = L^{\hp-1}_n$ and $K$ denote the respective
  lifts to $S$, we can see that $\hp_K - \hp > 0$ as follows.  First,
  we may assume that $c=2$, so that $n = d+2$,
  \[
  L = \llrr{x_0} + x_{1}^a \llrr{x_{1}, x_2, \ldots, x_{d}} \quad
  \text{ and } \quad K = x_{0}\llrr{x_{0}, x_{1}, \ldots, x_{d}} +
  \llrr{x_{1}^a}.
  \]
  Next, for large $j$, $L_j$ is spanned by all degree $j$ monomials
  divisible by $x_0$ or $x_1^a$ except those in $x_1^a \llrr{x_{d+1},
    x_{d+2}}^{j-a}$, while $K_j$ is spanned by all degree $j$
  monomials divisible by $x_0$ or $x_1^a$ except those in $x_0
  \llrr{x_{d+1}, x_{d+2}}^{j-1}$.  Specifically, for $j \ge a+2$ and
  $0 \le i \le 1$, let $N_{j}$ denote the number of monomials of $S_j$
  divisible by $x_0$ or $x_1^a$.  Then
  \[
  \dim_{\kk} L_j - \dim_{\kk} K_j = N_{j} - j+a-1 - (N_{j} - j) = a-1
  \ge 2,
  \]
  so $\hp_K(j) \ge \hp(j)-1+2$ as desired.  This implies
  $\lift(\nabla(I)) = L$ and by Pardue's criterion we find again that
  $I$ is an expansion of $L$.  (The same dimension-count shows that
  $K$ is a third saturated strongly stable ideal with Hilbert
  polynomial $\hp$ when $a = 2$; see Remark \ref{rmk:3borels}.)

  \vspace{0.25em}

  \paragraph{(ii)(c)}
  Here $b = (d, d, \dotsc, d, b_{a+1}, 0)$ with $d > b_{a+1} \ge 2$.
  Lemma~\ref{lem:HpolyBorel} and Theorem~\ref{thm:SSSunique} imply
  $\nabla(I) = L^{\nabla(\hp)}_{n-1}$ and setting $L := L^{\hp -1}_n =
  \mathfrak{m}_{c-2} + x_{c-1}^a \langle x_{c-1}, x_c, \dotsc,
  x_{n'-1} \rangle$, where $n' = c + d - b_{a+1}$ is as in
  Lemma~\ref{lem:twoborels}, we see that $I_{a+2} \subset L_{a+2}$ has
  codimension $1$.  The elements
  \[
  x_0 x_n^{a+1}, x_1 x_n^{a+1}, \dotsc, x_{c-2} x_n^{a+1},
  x_{c-1}^{a+1} x_n, x_{c-1}^a x_c x_n, \dotsc, x_{c-1}^a x_{n'-1} x_n
  \]
  cannot all belong to $I$, so Pardue's criterion implies $x_{c-2}
  x_n^{a+1} \notin I$ or $x_{c-1}^a x_{n'-1} x_n \notin I$.  The first
  option implies $I$ is the expansion of $L$ at $x_{c-2}$.  The second
  implies that $I$ is lexicographic.
\end{proof}

\begin{proof}[Proof of Theorem~\ref{thm:main1}]
  Theorem~\ref{thm:main2} proves the claim.
\end{proof}

\begin{remark}
  The formation of $J^0$ from $J$ in the proof of (i)(b) can be framed
  as follows: $G(J)$ is taken as a set of \emph{Borel generators} of
  the strongly stable ideal $J^0$.  Certain properties of strongly
  stable ideals can be described directly from their minimal Borel
  generators; see \cite{Francisco--Mermin--Schweig--2011}.  The idea
  of forming the strongly stable closure of a nonstandard Borel-fixed
  ideal is touched on briefly in \cite[Chapter VI, \S
    3]{Pardue--1994}; we are unaware of other instances of its use.
\end{remark}
  
\begin{remark}
  \label{rmk:3borels}
  Continuing with our method of adapting Reeves' algorithm to the
  Borel-fixed case, one can understand the Borel-fixed points on
  various Hilbert schemes with more than two Borel-fixed points in
  arbitrary characteristic.  For example, one can prove that the
  Hilbert scheme $\hilb^{\hp}(\PP^n)$ has exactly three Borel-fixed
  points, when $\hp$ has Gotzmann partition $(d, d, 1, 0)$ with $d >
  1$; cf.\ \cite[\S 4]{Ramkumar--2019}.  As noted in the proof of
  (i)(a'), nonstandard Borel-fixed ideals make their first appearance
  in classifying Hilbert schemes with three Borel-fixed points.
\end{remark}

\begin{remark}
  \label{rmk:Bertone}
  Monomial ideals that satisfy the condition of Lemma
  \ref{lem:HpolyBorel}(i) are known as \emph{\bfseries quasi-stable}
  or \emph{\bfseries of Borel type}.  These are also defined through a
  combinatorial criterion: for all monomials $m \in I$ and all $0 \le
  i < \max m$, there exists $s \ge 0$ such that $x_{\max m}^{-1} x_i^s
  m \in I$.  Stable ideals (those satisfying the same criterion but
  with $s = 1$) and $p$-Borel ideals are quasi-stable; see
  \cite[Proposition 4.2.9]{Herzog--Hibi--2011}.  Quasi-stable ideals
  are distinguished by the existence of \emph{Pommaret bases}
  \cite[Proposition 4.3]{Seiler--2009}, which are a special type of
  Stanley decomposition---the unique factorizations of the monomials
  of $J^0$ used in (i)(b) are an instance of this.  The existence of
  Pommaret bases allows Reeves' algorithm to be adapted to generate
  quasi-stable and $p$-Borel ideals \cite{Bertone--2015}.  Our proof
  of Theorem \ref{thm:main1} does not require the full machinery of
  this algorithm, however our approach independently identifies key
  properties applied in the algorithm.  (We thank C. Bertone for
  pointing out the connection to her work.)
\end{remark}

\begin{example}
  Let $n=3$. The ideal $I := \llrr{x_0^2, x_1} \subset R$ is
  quasi-stable with Hilbert polynomial $2$, while $\lift(I) =
  \llrr{x_0^2, x_1} \subset S$ is quasi-stable with Hilbert polynomial
  $2t+1$ (easily checked in e.g.\ \emph{Macaulay2}
  \cite{Grayson--Stillman}).  Neither ideal is Borel-fixed, thus
  Theorem \ref{thm:SSSunique} does not specify where there are unique
  quasi-stable ideals.  Moreover, the ideal $\llrr{x_0^2, x_0x_1,
    x_1^2, x_1x_2} \subset S$ is stable with Hilbert polynomial
  $2t+2$; again it is not Borel-fixed, so Theorem \ref{thm:main1} does
  not specify where exactly two quasi-stable ideals exist.
\end{example}

\begin{example}
  Let $f := x_1^2 + x_1x_2 + x_2^2 \in \FF_2[x_0, x_1, x_2, x_3]$.
  The ideal $\llrr{x_0, x_1f, x_2f}$ is Borel-fixed with Hilbert
  polynomial $2t+2$; in other words, Theorem \ref{thm:main1} does not
  hold over finite fields.
\end{example}

\begin{remark}
  \label{rmk:SkSm}
  Concurrent work by Skjelnes--Smith (initially inspired by our
  geography described in \S\ref{ch:geo}) has resulted in a full
  classification of smooth Hilbert schemes
  \cite{Skjelnes--Smith--2020}; their cases (4), (5), (6) correspond
  to cases (i)(c), (i)(b), (i)(a) here.  In fact, our results show
  that the difficulties associated with nonstandard Borel-fixed ideals
  in positive characteristic are avoided by the classification.
\end{remark}
  

\bibliography{twoborels}{} \bibliographystyle{amsalpha}


\end{document}